\documentclass[final]{elsarticle}

\usepackage[margin=2cm]{geometry}
\usepackage[polish,english]{babel}  
\usepackage[cp1250]{inputenc}
\usepackage{color}

\usepackage{hyperref}
  \hypersetup{
     	bookmarksopen=true
		,colorlinks = true
		, urlcolor = blue
  }

\usepackage{tikz}
\usepackage{subfig}

\usepackage{synttree}		
\usepackage{amsmath}

\newcommand{\pot}{u}
\newcommand{\der}{q}

\newcommand{\nonss}[1]{{#1}^*}

\newcommand{\ser}{S}

\DeclareMathAlphabet{\mathpzc}{OT1}{pzc}{m}{it}
\newcommand{\w}{\mathpzc{w}}

\newcommand{\bc}[1]{\bar{#1}}

\newcommand{\bound}{\Gamma}

\newcommand{\dif}{\mathrm{d}}

\newcommand{\vct}{\mathbf}

\definecolor{piggy}{RGB}{255,204,204}
\definecolor{perpel}{RGB}{204,204,255}
\definecolor{muddygreen}{RGB}{204,255,204}

\newenvironment{caparr}
{\begin{tabular}}
{\end{tabular}}

\begin{document}
\begin{frontmatter}

\title{The scalability of the matrices in direct Trefftz method in 2D Laplace problem}

\author{M. Borkowski} 

\ead{marbor\ @\ prz.edu.pl}

\address{Faculty of Electrical and Computer Engineering, Rzeszow University of Technology, Rzeszów, Poland}

\begin{abstract}
This paper presents an interesting property of the matrices that may be obtained with the use of direct Trefftz method. 
It is proved analytically for 2D Laplace problem that values of the elements of matrices describing the capacitance of two scaled domains are inversely proportional to the scalability factor. 
As an example of the application the capacitance extraction problem is chosen. Concise description of the algorithm in which the scalability property can be utilized is given. Furthermore some numerical results of the algorithm are presented.
\end{abstract}

\begin{keyword}
direct Trefftz method \sep 2D potential problem \sep scalable stiffness matrix

\end{keyword}

\end{frontmatter}


\section{Introduction}
While the Trefftz approach dates back to 1926 \cite{a123}, in literature it still mainly occurs in its indirect versions. 
The original Trefftz idea comes down to solving boundary value problem by the superposition of homogeneous solutions of the governing equation.
With the use of weighted residual method and employing different weighting functions one can obtain different versions of indirect Trefftz method, i.e. collocation, Galerkin or least squares version \cite{a6,b41,a405,a446}.
Systems of equations constituted by these methods are then solved for unknown parameters, so that solution approximated by the series of linearly independent functions satisfies the boundary conditions.

Indirect Trefftz method is a base for another, relatively new approach i.e., 
Hybrid Trefftz Method \cite{a22,a91,a29,b187,a360}, which can be regarded as the combination of FEM-like discretization with Trefftz-type approximation of the solution inside each finite element. 
Thanks to that, the governing equation is satisfied everywhere inside the domain of interest except for element boundaries. 
The interelement continuity is achieved by the additional approximation of the physical quantities on the element boundaries. 

While indirect Trefftz method is based on original variational formulation of the classical problem, its direct counterpart (dTM) is derived from the inverse one, just like popular direct BEM (dBEM) \cite{a233}.
The difference is while in dBEM the fundamental solutions of the equation are taken as a weighting functions, in dTM set of homogeneous solutions, called T-complete functions are used for that purpose.
Such a choice of weighting function results in a regular boundary integral equation. 

One of the merits of dTM is that unknowns are the actual physical variables of the problem and can be used to interpolate the physical quantity on the boundary.
On the other hand, the values in the interior of the domain must be obtained in terms of other procedures e.g. they can be calculated as a sum of single and double layer potentials. 
The final matrix obtained with the use of the method is full, nonsymmetric, non-diagonal dominant and very often ill-conditioned. 

Ill-conditioning of the matrices obtained by dTM may be the cause that the method is not so popular as other boundary procedures.
Nonetheless, it has been already applied for
potential \cite{a264,a20,a8,a404}, plane elasticity \cite{a265}, Kirchhoff plate bending \cite{a266} or free-vibration problems \cite{a137}.
There are also some works that concentrate on studying conditioning of the matrices of the method \cite{a20,a69,a19,a310,b213}.

The dTM still needs a further study as it has some its characteristic merits and in this paper we would like to present one that has not been noticed until now.
As it will be shown in next sections, the values of the matrices obtained in dTM for pure Dirichlet problems for two contours, where one of them is scaled with the factor $\lambda$, are inversely proportional to this factor.
The property is specific for direct Trefftz method and it can not be found in other boundary methods. 

The aim of this paper is to present the idea of the scalable stiffness matrices in dTM and to show its possible application. We also give some examples that confirm theoretical considerations.


 \section{Direct Trefftz method \label{Sec:Direct_Trefftz_Method}} 
 
 For simplicity, we consider the
 potential problem in the interior domain~$\Omega$ with 
 pure Dirichlet boundary conditions prescribed on piecewise regular boundary~$\bound$:
\begin{eqnarray} 
 \nabla^2 \pot = \frac{\partial^2 \pot}{\partial x^2} + \frac{\partial^2 \pot}{\partial y^2} = 0 & & \textnormal{in $\Omega$}
\label{Eq1:Lap-eq} \\
 \pot = \bc{\pot} & & \textnormal{on $\bound$} 
\label{Eq2:pot-BC} 
\end{eqnarray}
By making use of the weighted residual method one can obtain inverse variational formulation of the problem which is the starting point for direct boundary methods in general and dTM in particular:
\begin{equation}
	\int_\Omega  \pot \nabla^2\w \dif \Omega - \int_\bound \pot \frac{\partial \mathpzc{w}}{\partial n} \dif \bound + \int_\bound  \der \mathpzc{w} \dif \bound  = 0
\label{Eq:IWRF}
\end{equation} 
where $\w$ is a suitable weighting function, $\pot$ is prescribed by Eq.~(\ref{Eq2:pot-BC}) and $\der=\frac{\partial \pot}{\partial n}$ is unknown function on boundary.
In direct boundary methods $\pot$ and $\der$ can be interpreted physically, respectively as the potential and the electric flux on $\bound$.

The idea of dTM is to employ non-singular T-functions (T-complete functions) \cite{a353} for specific problem as the weighting function $\w$ in~Eq.~(\ref{Eq:IWRF}). For 2D potential problem in the interior domain they take the form
\begin{equation}
\nonss{\pot} \in  \left\{ 1, \rho^m \cos{m\theta}, \rho^m \sin{m\theta} \right\}
\label{Eq:TM-T-complete-funs}
\end{equation}
where $m=1,2,3\ldots$, and $(\rho,\theta)$ is the polar coordinate.
Since 
\begin{equation}
\nabla^2\nonss{\pot} =  0
\end{equation}
thus one can obtain regular boundary integral equation
\begin{equation}
	 \int_\bound \pot \nonss{\der}  \dif \bound = \int_\bound  \der \nonss{\pot} \dif \bound 
\label{Eq:regular_BIE}
\end{equation} 
where 
\begin{equation}
\nonss{\der} = \frac{\partial \nonss{\pot}}{\partial n} = \frac{\partial \nonss{\pot}}{\partial \rho} \frac{\partial \rho}{\partial n} +
\frac{\partial \nonss{\pot}}{\partial \theta} \frac{\partial \theta}{\partial n}
\label{Eq:T-fuN_cormal_deriv}
\end{equation}

Next, the procedure of finding potential $\pot$ and flux $\der$ requires boundary discretization and approximation of unknown functions with suitable interpolation polynomials on each element.
Let us assume in this paper that $\bound$ is divided into $N_j$ elements
\begin{equation}
\bound =  \cup_{j=1}^{N_j} \bound_j
\end{equation}
and the geometry of $j$-th element is defined with polynomials of arbitrary degree~$N_k$ in terms of local coordinate $\xi \in [-1, 1]$
\begin{equation}
	x(\xi)|_{\bound_j} = \sum_{k=1}^{N_k} x_k \ser_{rk}(\xi) 
	\textnormal{\ \ \ and\ \ \ }
	y(\xi)|_{\bound_j} = \sum_{k=1}^{N_k} y_k \ser_{rk}(\xi) \label{Eq:dBM-discretisation}
\end{equation}
where $\ser_{rk}$ is $k$-th interpolation polynomial of $j$-th boundary element.
The geometry can be expressed equivalently in polar coordinates since
\begin{equation}
	\rho = \sqrt{x^2+y^2} 	
	\textnormal{\ \ \ and\ \ \ }
	\theta = \arctan \frac{y}{x}  
	\label{Eq:cartesian_2_polar}
\end{equation}
The degree $N_\nu$ of polynomials that are used for approximation of $\pot$ and $\der$ on each element is also arbitrary. Thus
\begin{equation}
	\pot(\xi)|_{\bound_j} = \sum_{\nu=1}^{N_\nu} \pot_\nu^{(j)} \ser_{a\nu}(\xi) 	
	\textnormal{\ \ \ and\ \ \ }
	\der(\xi)|_{\bound_j} = \sum_{\nu=1}^{N_{\nu}} \der_\nu^{(j)} 
	\ser_{a\nu}(\xi) \label{Eq:boundary-function-approx}
\end{equation}
and $\pot_\nu^{(j)}$, $\der_\nu^{(j)}$ are respectively potential and flux values in $\nu$-th node on $j$-th element.
Hence, Eq.~(\ref{Eq:regular_BIE}) constitutes the system of linear equations where $i$-th equation is given by
\begin{equation}
\sum_{j=1}^{N_j} \sum_{\nu=1}^{N_\nu} \pot_\nu^{(j)}
 \int_{\bound_j}  \ser_{a\nu}(\xi) \nonss{\der}_i \dif \bound_j 
= \sum_{j=1}^{N_j} \sum_{\nu=1}^{N_\nu} \der_\nu^{(j)}
\int_{\bound_j} \ser_{a\nu}(\xi) \nonss{\pot}_i \dif \bound_j
\end{equation}
The system of equations may be written in matrix form as
\begin{equation}
\vct{ Hu=Gq }
\label{Eq:Hu=Gq}
\end{equation}
where $\vct u = [u_1, u_2, \ldots, u_{N_t}]^T$, $\vct  q = [q_1, q_2, \ldots, q_{N_t}]^T$ ($N_t = N_j \cdot N_k$ -- the total number of interpolation nodes on boundary), and elements of matrices $\vct H$ and $\vct G$ take the form
\begin{eqnarray} 
h_{ij}^{(k)} = \int_{\bound_j} \ser_{a\nu} \nonss{\der}_i \dif \bound_j \label{Eq:Hij} \\
g_{ij}^{(k)} = \int_{\bound_j} \ser_{a\nu} \nonss{\pot}_i \dif \bound_j \label{Eq:Gij}
\end{eqnarray}
Taking into consideration boundary conditions~(\ref{Eq2:pot-BC}), the solution of Eq.~(\ref{Eq:Hu=Gq}) is equal
\begin{equation}
\vct \der = \vct C \vct {\bc{\pot}}
\label{Eq:BCM-Equation}
\end{equation}
where
\begin{equation}
\vct C = \vct G^{-1} \vct H
\label{Eq:C-matrix}
\end{equation}
It is obvious that the number of weighting functions should be equal to~$N_t$ to compute inverse of~$\vct G$.

Obtained in such a way matrix $\vct C$ can be used to combine boundary method with FEM \cite{b21,b36,b55}. Some authors \cite{a404,a146,a148} call it \emph{boundary capacitance matrix} (BCM) to emphasize its resemblance to the matrices describing coupling capacitances of the systems of conductors. $\vct C$ elements express charge induced on boundary elements in terms of their potentials. 
However, it should be emphasized that this interpretation has its physical grounds only when dBEM is used for the calculation of $\vct C$. Then, the fundamental solutions and their normal derivatives employed as weighting functions correspond to the effect of the boundary load acting on its neighbourhood.
When T-complete solutions are employed in dTM $\vct C$ elements \emph{does not have} this physical interpretation.


\section{Calculation of coefficients of $\vct H$ and $\vct G$}

\begin{figure}
\centering
{ \input{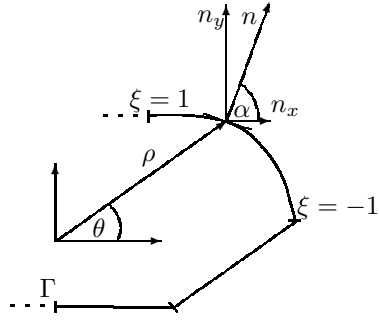}}
\caption{Boundary element parametrization}
\label{Fig:boundary_parametrization}
\end{figure}

Now let us take a closer look at the process of computing $\vct H$ and $\vct G$ to point out the property of $\vct C$ obtained in dTM. 

Curvilinear integrals in Eqs.~(\ref{Eq:Hij}) and (\ref{Eq:Gij}) should be reduced to definite integrals. Derivation of~(\ref{Eq:dBM-discretisation}) gives:
\begin{equation}
	\frac{\dif x}{\dif \xi}|_{\bound_j} = \sum_{k=1}^{N_k} x_k \frac{\dif \ser_{rk}}{\dif \xi} 
	\textnormal{\ \ \ and\ \ \ }
	\frac{\dif y}{\dif \xi}|_{\bound_j} = \sum_{k=1}^{N_k} y_k \frac{\dif \ser_{rk}}{\dif \xi} 
	\label{Eq:cartesian-coordinates-derivatives}
\end{equation}
which can be used to obtain the Jacobian of the transformation between local and global coordinates:
\begin{equation}
\dif \bound_j = \sqrt{\dif x^2 + \dif y^2} =J_j \dif \xi
\end{equation}
Taking into account (\ref{Eq:cartesian_2_polar}) one can express normal derivatives of $\rho, \theta$ from (\ref{Eq:T-fuN_cormal_deriv}) as
\begin{equation}
	\frac{\partial \rho}{\partial n} = 
	\frac{\partial \rho}{\partial x} \frac{\partial x}{\partial n} + 
	\frac{\partial \rho}{\partial y} \frac{\partial y}{\partial n}
	\textnormal{\ \ \ and\ \ \ }
	\frac{\partial \theta}{\partial n} = 
	\frac{\partial \theta}{\partial x} \frac{\partial x}{\partial n} + 
	\frac{\partial \theta}{\partial y} \frac{\partial y}{\partial n}	
	\label{Eq:cartesian_2_polar-derivatives}
\end{equation}
where
\begin{equation}
	\frac{\partial \rho}{\partial x} = \frac{x}{\rho}
	\textnormal{, \ \ \ }	
	\frac{\partial \rho}{\partial y} = \frac{y}{\rho}
	\textnormal{, \ \ \ }	
	\frac{\partial \theta}{\partial x} = -\frac{y}{\rho^2}
	\textnormal{\ \ \ and\ \ \ }	
	\frac{\partial \theta}{\partial y} = \frac{x}{\rho^2}
	\label{Eq:polar-cartesian-partials}
\end{equation}
Since
\begin{equation}
	\frac{\partial x}{\partial n} = \cos \alpha = \frac{\dif y}{\dif \bound_j} 
	\textnormal{\ \ \ and\ \ \ }
	\frac{\partial y}{\partial n} = -\sin \alpha = -\frac{\dif x}{\dif \bound_j} 
	\label{Eq:cartesiaN_cormal_derivatives}
\end{equation}
and taking into account polar to Cartesian coordinate transformation
\begin{equation}
	x = \rho \cos \theta
	\textnormal{\ \ \ and\ \ \ }
	y = \rho \sin \theta
	\label{Eq:polar_2_cartesian}
\end{equation}
the components of the outward normal at a boundary element can be given in terms of polar coordinates and a corner angle $\alpha$ (Fig.~\ref{Fig:boundary_parametrization}) as
\begin{eqnarray} 
& & \frac{\partial \rho }{\partial n} = 
\cos \theta \cos \alpha - \sin\theta \sin \alpha =
\cos (\theta + \alpha)
\label{Eq:dr_dn}
\\
& & \frac{\partial \theta }{\partial n} = 
\frac{1}{\rho} \left( -\sin \theta \cos \alpha - \cos\theta \sin \alpha \right) = -\frac{1}{\rho} \sin(\theta + \alpha)
\label{Eq:dt_dn}
\end{eqnarray}
Let us observe, that taking into account above formulae Eq.~(\ref{Eq:T-fuN_cormal_deriv}) 
and trigonometric 
angle sum and difference identities,
$\nonss{\der}$  can also be expressed as a set of functions analogously to set~(\ref{Eq:TM-T-complete-funs})
\begin{equation}
\nonss{\der} \in  \left\{ 0, 
m\rho^{m-1} \cos \left( \left( m-1 \right) \theta -\alpha \right),
m\rho^{m-1} \sin \left( \left( m-1 \right) \theta -\alpha \right)
\right\}
\label{Eq:TM-T-complete-deriv-funs}
\end{equation}
For the sake of brevity of further formulae let us express $\nonss{\pot}$, $\nonss{\der}$ respectively as:
\begin{equation}
	\nonss{\pot} = R^u(\rho) T^u(\theta)
	\textnormal{\ \ \ and\ \ \ }
	\nonss{\der} = R^q(\rho) T^q(\theta,\alpha)
\end{equation}
where again
\begin{alignat}{2} 
& R^u \in  \left\{ 1, \rho^m, \rho^m \right\} 
& \textnormal{\ \ and\ \ } 
& T^u \in  \left\{ 1, \cos{m\theta}, \sin{m\theta}  \right\}
\label{Eq:R_u-T_u-sets} \\
& R^q \in  \left\{0, \rho^{m-1}, \rho^{m-1}  \right\} 
& \textnormal{\ \ and\ \ } 
& T^q \in  \left\{ 
0, 
m \cos \left( \left( m-1 \right) \theta -\alpha \right), 
m \sin \left( \left( m-1 \right) \theta -\alpha \right)  
\right\}
\label{Eq:R_q-T_q-sets}
\end{alignat}

Hence, one can express Eqs.~(\ref{Eq:Hij}) and (\ref{Eq:Gij}) as:
\begin{eqnarray} 
h_{ij}^{(k)} = \int_{-1}^{1} \ser_{a\nu}  
R^q_i T^q_i
J_j \dif \xi \label{Eq:Hij_local} \\
g_{ij}^{(k)} = \int_{-1}^{1} \ser_{a\nu} 
R^u_i T^u_i
J_j \dif \xi \label{Eq:Gij_local}
\end{eqnarray}

\section{Scalability of the matrices of the scaled boundaries}
\begin{figure}
\centering
{ \input{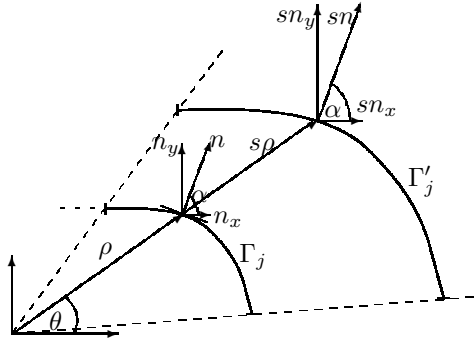}}
\caption{The scaling of the boundary element}
\label{Fig:boundary_scaling}
\end{figure}

Now let us consider a uniform scaling of the domain $\Omega$ and boundary $\Gamma$  
where the nonzero scaling factor is equal $s$:
\begin{equation}
\Omega' \rightarrow s \Omega
\end{equation}
which means that each point $(x_k,y_k)$ from $\Omega$ and $\bound$ is transformed into $(sx_k,sy_k)$ (Fig.~\ref{Fig:boundary_scaling}). 
In a similar way as presented above one can derive elements of matrices $\vct H'$ and $\vct G'$ of the scaled boundary to be:
\begin{eqnarray} 
h_{ij}^{'(k)} 
= \int_{-1}^{1} \ser_{a\nu}  R^q_i(s\rho) T^q_i s J_j \dif \xi 
\label{Eq:scaled_Hij_local} \\
g_{ij}^{'(k)} 
= \int_{-1}^{1} \ser_{a\nu}  R^u_i(s\rho) T^u_i s J_j \dif \xi  
\label{Eq:scaled_Gij_local}
\end{eqnarray}
Let us observe that taking into account the form of $R$ functions the following formulae hold true:
\begin{equation}
	R^u(s\rho) = R^u(s)R^u(\rho)
	\textnormal{\ \ \ and\ \ \ }
	R^q(s\rho) = R^q(s)R^q(\rho)
\end{equation}
and Eqs.~(\ref{Eq:scaled_Hij_local}) and (\ref{Eq:scaled_Gij_local}) can be rewritten as
\begin{eqnarray} 
h_{ij}^{'(k)} 
= sR^q_i(s) \int_{-1}^{1} \ser_{a\nu} R^q_i(\rho) T^q_i J_j \dif \xi
= sR^q_i(s) h_{ij}^{(k)}
\label{Eq:scaled_Hij_local2} \\
g_{ij}^{'(k)} 
= sR^u_i(s) \int_{-1}^{1} \ser_{a\nu}(\xi)  R^u_i(\rho) T^u_i J_j \dif \xi 
= sR^u_i(s) g_{ij}^{(k)} 
\label{Eq:scaled_Gij_local2}
\end{eqnarray}
Thus one can express $\vct H'$ and $\vct G'$ as
\begin{equation}
\vct H' = \left[ 
\begin{array}{cccc}
	   sR^q_1(s) h_{11}^{(1)}    &    sR^q_1(s) h_{11}^{(2)}    & \ldots &    sR^q_1(s) h_{1N_j}^{(N_k)}    \\
	   sR^q_2(s) h_{21}^{(1)}    &    sR^q_2(s) h_{21}^{(2)}    & \ldots &    sR^q_2(s) h_{2N_j}^{(N_k)}    \\
	          \ldots           &           \ldots           & \ldots &             \ldots             \\
	sR^q_{N_t}(s) h_{N_t1}^{(1)} & sR^q_{N_t}(s) h_{N_t1}^{(2)} & \ldots & sR^q_{N_t}(s) h_{N_tN_j}^{(N_k)}
\end{array}%
\right] 
\end{equation}

\begin{equation}
\vct G' = \left[ 
\begin{array}{cccc}
	   sR^u_1(s) g_{11}^{(1)}    &    sR^u_1(s) g_{11}^{(2)}    & \ldots &    sR^u_1(s) g_{1N_j}^{(N_k)}    \\
	   sR^u_2(s) g_{21}^{(1)}    &    sR^u_2(s) g_{21}^{(2)}    & \ldots &    sR^u_2(s) g_{2N_j}^{(N_k)}    \\
	           \ldots            &            \ldots            & \ldots &              \ldots              \\
	sR^u_{N_t}(s) g_{N_t1}^{(1)} & sR^u_{N_t}(s) g_{N_t1}^{(2)} & \ldots & sR^u_{N_t}(s) g_{N_tN_j}^{(N_k)}
\end{array}%
\right] 
\end{equation}
Both matrices $\vct H'$ and $\vct G'$ can be decomposed into two parts and expressed 
as products of some diagonal matrices, which depend on a scale factor $s$, and respectively $\vct H$ and $\vct G$ matrices of the original, unscaled boundary $\bound$:
\begin{equation}
	\vct H' = \vct H_s \vct H 
	\textnormal{\ \ \ and\ \ \ }
	\vct G' = \vct G_s \vct G \label{Eq:scaled_H_and_G_matrices}
\end{equation}
where
\begin{equation}
\vct H_s = \left[ 
\begin{array}{cccc}
	sR^q_1(s) &     0     & \ldots &       0         \\
	    0     & sR^q_2(s) & \ldots &       0         \\
	 \ldots   &  \ldots   & \ldots &    \ldots       \\
	    0     &     0     & \ldots & sR^q_{N_t}(s) 
\end{array}%
\right] 
	\textnormal{\ \ \ and\ \ \ }
\vct G_s = \left[ 
\begin{array}{cccc}
	sR^u_1(s) &     0     & \ldots &       0         \\
	    0     & sR^u_2(s) & \ldots &       0         \\
	 \ldots   &  \ldots   & \ldots &    \ldots       \\
	    0     &     0     & \ldots & sR^u_{N_t}(s) 
\end{array}%
\right] 
\end{equation}
Now, taking into account that $\vct{(AB)}^{-1} = \vct B^{-1} \vct A^{-1}$, $\vct C'$ can be calculated as
\begin{equation}
\vct C' =  
\vct G'^{-1} \vct H' = 
\vct G^{-1} \vct G_s^{-1} \vct H_s \vct H 
\label{Eq:C-for-scaled-matrix}
\end{equation}
thus
\begin{equation}
\vct C' =  \vct G^{-1} \vct S \vct H 
\end{equation}
where $\vct S = \vct G_s^{-1} \vct H_s$.
From Eqs. (\ref{Eq:R_u-T_u-sets}) and (\ref{Eq:R_q-T_q-sets}) 
one can conclude that $i$-th $R$ functions are correspondingly functions taken from the sets:
$$
\begin{array}{cccll}
	R^u_i(s) \in \left\{ \right. & 1,  s,  s, & s^2, & s^2, & s^3, s^3, s^4, \ldots \left. \right\}  \\
	R^q_i(s) \in \left\{ \right. & 0,  1,  1, & s,   & s,   & s^2, s^2, s^3, \ldots \left. \right\}  
\end{array}
$$
Thus it can be seen that:
$$
\frac{R^q_1(s)}{R^u_1(s)} = 0
	\textnormal{\ \ \ and\ \ \ }
\frac{R^q_i(s)}{R^u_i(s)} = \frac{1}{s} 
	\textnormal{ for } i = 2,3, \ldots
$$
and
\begin{equation}
\vct S 
 = \left[ 
 \begin{array}{cccc}
 	\frac{1}{sR^u_1(s)} &     0     & \ldots &       0       \\
 	    0     & \frac{1}{sR^u_2(s)} & \ldots &       0       \\
 	 \ldots   &  \ldots   & \ldots &    \ldots     \\
 	    0     &     0     & \ldots & \frac{1}{sR^u_{N_t}(s)} 
 \end{array}%
 \right] 
 \left[ 
 \begin{array}{cccc}
 	sR^q_1(s) &    0    & \ldots &      0      \\
 	   0    & sR^q_2(s) & \ldots &      0      \\
 	\ldots  & \ldots  & \ldots &   \ldots    \\
 	   0    &    0    & \ldots & sR^q_{N_t}(s)
 \end{array}%
 \right] 
  = \frac{1}{s} \vct I_p
\end{equation}
and $\vct I_p$ is the identity matrix or diagonal matrix made from identity matrix where $p$-th diagonal element can be set to $0$, which reflects the fact that one of the weighting functions can be the first one from the set~(\ref{Eq:TM-T-complete-funs}).
Having in mind that in such case $p$-th row of $\vct H$ is a zero vector one can write
\begin{equation}
\vct C' = \vct G^{-1} \frac{1}{s} \vct I_p \vct H = \frac{1}{s} \vct G^{-1} \vct H = \frac{1}{s} \vct C
\end{equation}
Hence, the elements of $\vct C$ matrices for two given domains, from which one is an isotropic transformation of another, are inversely proportional to the transformation factor.

Some calculation examples and the possible numerical application of this feature are presented in the next section.


\section{Example problems}

\subsection{Exact formulae for simple geometries} \label{Sec:Exact-formulae}
Consider some simple geometries of  constant boundary elements that can be generated for rectangular domain 
$\Omega = \left\{(x,y), a < x < a+w, b< y < b+h \right\}$, where the symbols are depicted in {Fig.~\ref{Fig:Example-boundaries}. 

For the convenient calculations, T-functions from the set given by Eq.~(\ref{Eq:TM-T-complete-funs}) can be equivalently expressed in Cartesian coordinates as:
\begin{equation}
\nonss{\pot} \in  \left\{ 1, x, y, x^2-y^2, 2xy, x^3-3xy^2, \ldots \right\}
\label{TM-T-complete-funsXY}
\end{equation}

\begin{figure}[b]
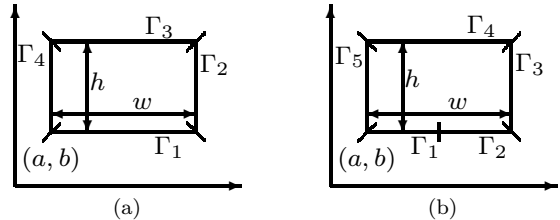

\centering
\subfloat[]{
\scalebox{1.2}
{ \footnotesize
\input{TexFig/element1.texfig}}
}
\qquad
\subfloat[]{
\scalebox{1.2}
{ \footnotesize
\input{TexFig/element2.texfig}}
}
 \caption{Boundary discretization (a) for Examples 1 and 2; (b) for Example 3} 
\label{Fig:Example-boundaries}
\end{figure}

\subsubsection*{Example 1}
Let us consider the geometry presented in Fig.~\ref{Fig:Example-boundaries}(a). Taking the first four functions from~(\ref{TM-T-complete-funsXY}) as the weighting functions one can easily obtain
\begin{equation}
\vct{H} = 
\left[
\begin{array}{c c c c}
	 0  & 0       & 0        &  0       \\
	 0  & h       & 0        & -h       \\
	-w  & 0       & w        &  0       \\
	2bw & 2h(a+w) & -2w(b+h) & -2ah
\end{array}
\right]
\label{Eq:Ex1-H}
\end{equation}
and 
\begin{equation}
\vct{G} = 
 \left[
 \arraycolsep=3.0pt\def\arraystretch{1.4}
\begin{array}{cccc}
	         w          & h                       & w                       & h                   \\
	\frac{1}{2}w^2 +aw  & h(a+w)                  & \frac{1}{2}w^2 +aw  & ah                  \\
	        bw          & \frac{1}{2}h^2 +bh      & w(b+h)                  & \frac{1}{2}h^2 +bh  \\
g_{41}	 & g_{42} & g_{43} & g_{44}
\end{array}
\right]
\label{Eq:Ex1-G}
\end{equation}
where 
$g_{41} = w(a^2+aw-b^2+\frac{w^2}{3})$,
$g_{42} = h[(a+w)^2-b^2-bh-\frac{h^2}{3}]$,
$g_{43} = w[a^2+aw-(b+h)^2+\frac{w^2}{3}]$,
$g_{44} = h(a^2-b^2-bh-\frac{h^2}{3})$.

It can be seen that the elements of both $\vct{H}$ and $\vct{G}$ are expressed in terms of both domain position ($a$,$b$), which is equivalent to the position of the centre of eigenexpansion, and its size ($w$,$h$). On the other hand, $\vct{C}$ for this domain equals:
\begin{equation}
\vct{C} = \left[
\begin{array}{c c c c}
\frac{4h^2 + w^2}{h^3+hw^2} & \frac{-3h }{h^2+w^2} & \frac{2h^2-w^2}{h^3+hw^2} & \frac{-3h }{h^2+w^2} \\
\frac{ -3w }{h^2+w^2} & \frac{h^2+4w^2}{h^2w+w^3} & \frac{-3w }{h^2+w^2} & \frac{2w^2 - h^2}{h^2w+w^3} \\
\frac{2h^2 - w^2}{h^3+hw^2} & \frac{      -3h }{h^2+w^2} & \frac{4h^2 + w^2}{h^3+hw^2} & \frac{      -3h }{h^2+w^2} \\
\frac{      -3w}{h^2+w^2} & \frac{2w^2 - h^2}{h^2w+w^3}  & \frac{      -3w }{h^2+w^2} & \frac{h^2+4w^2}{h^2w+w^3}
\end{array}
\right]
\label{Eq:Ex1-C}
\end{equation}  
and its elements \emph{do not} depend on the position of the domain on the plane. 

\subsubsection*{Example 2}
The same $\vct C$ matrix can be obtained when other combinations of functions from Eq.~(\ref{TM-T-complete-funsXY}) are taken to be weighting functions. As an example let us take four consecutive functions beginning from the second one i.e. $\w \in  \left\{x, y, x^2-y^2, 2xy \right\}$. For the sake of conciseness $\vct{H}$ and $\vct{G}$ matrices are not presented for this example.

Though $\vct{H}$ and $\vct{G}$ for this example are different from Eqs.~(\ref{Eq:Ex1-H}) and (\ref{Eq:Ex1-G}), $\vct C$ is exactly the same as in Eq.~(\ref{Eq:Ex1-C}).

\subsubsection*{Example 3}
In this example let us consider a bit more complicated boundary geometry presented in~Fig.~\ref{Fig:Example-boundaries}(b). Analogously to both above examples one can obtain $\vct H = \vct H(a,b,w,h)$, $\vct G = \vct G(a,b,w,h)$ and 
\begin{equation}
\vct C = \vct C(w,h) = 
\left[
\begin{array}{c c c c c}
\frac{6h^2+3w^2}{2h^3+2hw^2} & \frac{2h^2-w^2}{2h^3+2hw^2} & \frac{-3h }{h^2+w^2} & \frac{2h^2-w^2}{h^3+hw^2} & \frac{-3h }{h^2+w^2} \\
\frac{2h^2-w^2}{2h^3+2hw^2} & \frac{6h^2+3w^2}{2h^3+2hw^2} & \frac{-3h }{h^2+w^2} & \frac{2h^2-w^2}{h^3+hw^2} & \frac{-3h }{h^2+w^2} \\
\frac{2w^3-wh^2}{4h^2(h^2+w^2)} & \frac{5wh^2+ 2w^3}{-4h^2(h^2+w^2)} & \frac{h^2+4w^2}{h^2w+w^3} & \frac{-3w }{h^2+w^2} & \frac{2w^2 - h^2}{h^2w+w^3} \\
\frac{2h^2 - w^2}{2(h^3+hw^2)} & \frac{2h^2 - w^2}{2(h^3+hw^2)} & \frac{      -3h }{h^2+w^2} & \frac{4h^2 + w^2}{h^3+hw} & \frac{      -3h }{h^2+w^2} \\
\frac{5wh^2+ 2w^3}{-4h^2(h^2+w^2)} &  \frac{2w^3-wh^2}{4h^2(h^2+w^2)} & \frac{2w^2 - h^2}{h^2w+w^3}  & \frac{      -3w }{h^2+w^2} & \frac{h^2+4w^2}{h^2w+w^3}
\end{array}
\right]
\end{equation}

These simple computational examples show that $\vct C$ matrices obtained by dTM are not only scalable for uniformly scaled geometries but also they do not depend either on the position of the centre of eigenexpansion or the selection of weighting functions. Thus, $\vct C$ can indeed be called a boundary capacitance matrix as it stores the information about the discretization of the domain boundary.
It should be noted that though in direct standard or regular BEM $\vct C$ matrix is also independent to the position of problem domain, it is not scalable in the above mentioned sense.

\subsection{Application of scalable Trefftz matrices for capacitance extraction}

\begin{figure}[b]
\centering
{\input{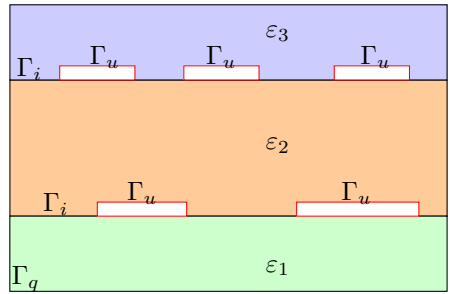}}
 \caption{Non-homogeneous multi-layer 2D planar transmission line as an example of system of $N_c$ conductors.} 
\label{Fig:transmissionLine}
\end{figure}

Fast and accurate estimation of parasitic capacitances is a very important aspect of modern integrated circuits designing \cite{a155,a160}. Due to the growing complexity of the geometry of systems of conductors new numerical approaches are sought to deal with the problem.
One of proposed solutions is the application of direct boundary methods.

The standard way of dealing with complexity of problem geometry in capacitance extraction is decomposition of the domain of the problem into elements that match the elements contained in library created during preprocessing stage. \cite{a148,a202}.


In recent paper we presented general algorithm that can be used with any direct boundary method and showed its application to the problem of non-homogeneous planar transmission lines parasitic capacitance calculation \cite{a404}. The geometry of example problem is presented in Fig.~\ref{Fig:transmissionLine}, where $\Gamma_u$ are boundaries of the conducting paths, $\Gamma_i$ -- interfaces between dielectric layers with different dielectric constants $\varepsilon_k$ ($k=1,2,3$), and $\Gamma_q$ is external boundary with homogenuous Neumann condition.

The general idea of the algorithm comes down to hierarchical combining $\vct C$ matrices (calculated by the means of appropriate direct boundary method) of adjacent subdomains which finally gives generalized capacitance matrix of the whole system. 
Here, we will focus only on the domain decomposition aspect of this algorithm to point out the possible application of the property of scaling Trefftz matrices.

The mesh generation method in \cite{a404} is hierarchical binary tree decomposition of each dielectric layer, where mesh density is higher near to the conducting paths. The example of such discretization and the corresponding binary tree are presented in Fig.~\ref{Fig:binary_tree_discretisation}.
The characteristic feature of hierarchical binary tree decomposition is that 
mesh consists of domains that most of can be classified into limited set of elements.

\begin{figure}[hbt]
\centering
\subfloat[]{
\scriptsize
\input{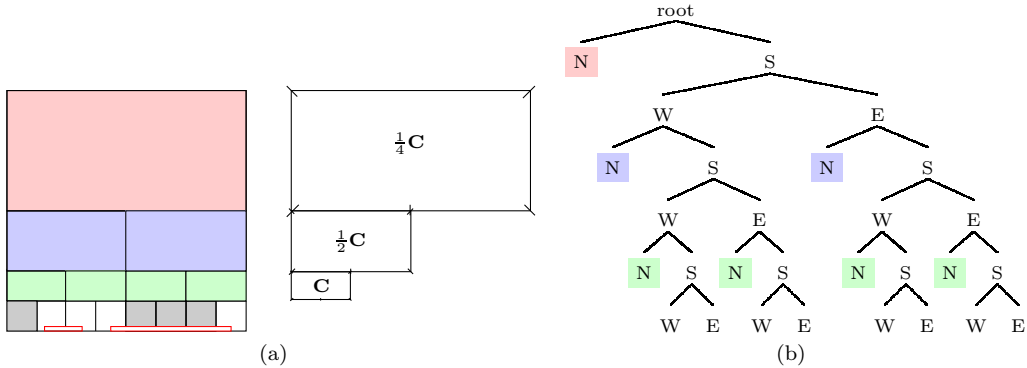}
}
\subfloat[]{
\scriptsize
\branchheight{.7cm}
\childsidesep{0em}
\childattachsep{2em}
\synttree[root
[ {\colorbox{piggy}{N}} ] [S [W~[ {\colorbox{perpel}{N}} ] [S [W [ {\colorbox{muddygreen}{N}} ] [S [W] [E]]] [E [{\colorbox{muddygreen}{N}}] [S [W] [E]]]]] [E [ {\colorbox{perpel}{N}} ] [S [W [{\colorbox{muddygreen}{N}}] [S [W] [E]]] [E [{\colorbox{muddygreen}{N}}] [S [W] [E]]]]]] 
]
}
\caption{(a) Domain decomposition of a single dielectric layer with two conducting paths at the bottom of the structure; (b) a corresponding binary tree (b).}
\label{Fig:binary_tree_discretisation}
\end{figure}

Now let us analyze the discretization depicted in Fig.~(\ref{Fig:binary_tree_discretisation})(a). It can be noticed that domain elements with the same boundary discretization may appear on different levels of the hierarchical binary tree (north subdomains coloured with three colours). Normally (when dBEM is applied as the PDE solver), they have to be calculated and stored separately as the BCMs are different for them.
The property of scalability of $\vct C$ matrix shown in previous sections allows to calculate and store $\vct C$ matrix only for the biggest (smallest) sub-domain of each type. The BCMs for the rest of the subdomains of this type can be obtained simply by multiplication of stored $\vct C$ matrix and the inverse of appropriate scaling geometrical factor.
Thus for the example shown in Fig.~\ref{Fig:binary_tree_discretisation}(a) instead of 3 calculations for north node domains of each size that appear on different binary tree levels only 1 calculation is needed. 
Let us also observe, that according to the examples presented in subsection \ref{Sec:Exact-formulae}, for very simple discretizations of the domain boundaries (such as those obtained in our problem) the analytical formulae for $\vct C$ matrices can be easily derived.

Of course, the deeper hierarchical tree/ more dense domain decomposition is, the better improvements in speed/memory requirements of the algorithm is gained.

\subsection{Numerical results} \label{Sec:numerical-results}
In presented examples the same assumptions as in \cite{a404} are made.
The quasi-transverse electromagnetic wave propagates along infinite conductive paths and the thickness of each path is negligible. Infinite paths length assumption implies that parasitic capacitance values are determined per unit length.
Dielectric layers are linear, homogeneous and isotropic. Results obtained with Linpar software \cite{b99} are taken as reference solutions ($\vct C_{ref}$). 

In below-mentioned examples the calculations are conducted with constant boundary elements.
The relative permittivities and geometry details are presented in the figures.
The heights of the dielectric layers, widths of the paths and distances between them are given in $mm$.
Conductors are numbered by the consecutive integers from left to right and from bottom to top layer.

The computations are performed on standard PC (Intel Core i5 2.3 GHz, 4 GB RAM) in Matlab 7.12.

\begin{figure}
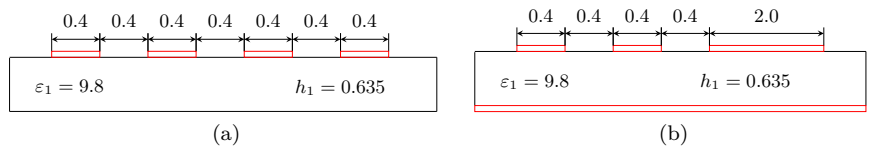

\centering
\subfloat[]{
{\small
\scalebox{0.8}{
\input{TexFig/p36_3.texfig} 
}}
}
\subfloat[]{
{\small
\scalebox{0.8}{
\input{TexFig/p36_5.texfig} 
}}
}
\caption{Examples of transmission lines from section \ref{Sec:numerical-results} }
\label{Fig:trans-lines-examples}
\end{figure}

\subsubsection*{Example 1}

In this example the first path on the left is taken as a ground. The generalized capacitance matrix is equal
$\vct C_{ref} = \left[
\begin{caparr}{rrr}
109.1   &  -46.0 & -10.5  \\
-46.0   &  109.1 & -52.6  \\
-10.5   & -52.6  & 70.6
\end{caparr} 
\right] pF/m$. 

The results are shown in Table \ref{Tab:Example1-results}. The consecutive columns correspond to higher and higher mesh density. The consecutive rows presents the total number of nodes on conductors
(total number of all nodes) -- $N_n$, generalized capacitance matrix
-- $\vct C_G$, root mean square error (RMSE) between obtained and reference solutions -- $E$, computation times with (new version of the algorithm) and without (original version) applying scalability of Trefftz elements respectively (in seconds) -- $t_n$ and $t_o$ correspondingly.

\begin{table}
$\begin{array}{|c|c|c|c|}
\hline
	   N_n     & 44 (597) & 84 (1037) & 168 (1845) \\ \hline 
	\vct C_G [pF/m] & 
\left[
\begin{caparr}{rrr}
107.1   &  -45.7 &  -8.7  \\
-45.7   &  107.1 & -52.7  \\
 -8.5   & -52.9  & 63.4
\end{caparr} 
\right]  & 
\left[
\begin{caparr}{rrr}
108.4   &  -46.3& -8.6 \\
-46.3   &  108.4 & -53.5  \\
 -8.5   & -53.6  & 64.1
\end{caparr} 
\right]
         &
\left[
\begin{caparr}{rrr}
109.4   &  -46.9 &  -8.6  \\
-46.9   &  109.4 & -53.5  \\
 -8.6   & -53.6  & 64.4
\end{caparr} 
\right]
\\ \hline  
	    E      & 7.4844   & 5.8700    & 5.6656 \\ \hline 
	 t_o [s]   & 0.7074   & 0.9592    & 1.3891 \\ \hline 
	 t_n [s]   & 0.6039   & 0.7996    & 1.2675 \\ \hline 
\end{array} 
$

\caption{Results obtained for the system presented in Fig.~\ref{Fig:trans-lines-examples}(a)}
\label{Tab:Example1-results}
\end{table}

\subsubsection*{Example 2}
In this example the ground plane is placed at the bottom of the structure. The reference solution is given by
$\vct  C_{ref} = \left[
\begin{caparr}{rrr}
143.3   &  -23.7 &   -1.8  \\
-23.7   &  148.1 &  -25.4  \\
 -1.8   &  -25.4 &  370.1
\end{caparr} 
\right] pF/m$.

The results are given in Table \ref{Tab:Example2-results} and the symbols used in the table have the same meaning as in the previous example.

\medskip

The solutions obtained with the newer, optimized version of the algorithm are the same as those obtained with the original version ($\vct C_G$). The advantage of the application of scalable Trefftz elements is a noticeable shortening of time consumption needed for the calculations -- $t_n$ is about 10\% lesser than $t_o$.

\begin{table}
$
\begin{array}{|c|c|c|c|}
\hline
	   N_n     & 332 (2069) & 663 (4073) & 1324 (8052) \\ \hline 
	\vct C_G [pF/m] & 
\left[
\begin{caparr}{rrr}
140.3   &  -22.6 &  -2.0  \\
-23.5   &  145.9 & -26.9  \\
 -1.8   & -34.3  & 375.2
\end{caparr} 
\right]  & 
\left[
\begin{caparr}{rrr}
141.5   &  -23.6 &  -1.5  \\
-23.6   &  146.4 & -25.4  \\
 -1.5   & -25.3  & 365.6
\end{caparr} 
\right]
         &
\left[
\begin{caparr}{rrr}
142.4   &  -24.1 &  -1.5  \\
-24.1   &  147.5 & -25.6  \\
 -1.5   & -25.6  & 366.1
\end{caparr} 
\right]
\\ \hline  
	    E     & 13.4844   & 2.9544    & 1.8844 \\ \hline 
	 t_o [s]   & 1.2060   & 1.9464   & 3.7753 \\ \hline 
	 t_n [s]   & 1.0696   & 1.7499    & 3.4774 \\ \hline 
\end{array} 
$
\caption{Results obtained for the system presented in Fig.~\ref{Fig:trans-lines-examples}(b)}
\label{Tab:Example2-results}
\end{table}

\section{Concluding remarks}
This paper presents the scalability of Trefftz elements derived in direct Trefftz method. This is done by introducing boundary capacitance matrix, which in fact, is a matrix obtained 
for pure Dirichlet problem. 

The scalability of boundary capacitance matrices is proved theoretically for 2D Laplace interior problem. Boundary capacitance matrices obtained for two scaled domains with the same boundary discretization are inversely proportional to the scalability parameter.

Next, some computational examples are given that confirm this proof. 
It seems from them, that above-mentioned scalability is not only independent of the combination of T-functions used as the weighting functions, but also it is independent of the location of centre of eigenexpansion, though these properties have not been proved yet.
It should be emphasized that scalability property is distinctive only for direct Trefftz method and other direct boundary methods (standard/regular BEM) does not have it.

Moreover, the possible application of the scalability in numerical calculations is demonstrated on the example of parasitic capacitance extraction of 2D planar structures. 
Presented examples show that the scalability of direct Trefftz elements can be utilized to speed up numerical calculations.
The further application of the scalable elements can be a combination with standard finite elements.

\section*{\footnotesize Acknowledgements}
\footnotesize
Numerical experiments were conducted with the use of MATLAB application, purchased during realisation of Project No. UDA-RPPK.01.03.00-18-003/10-00 "Construction, expansion and modernisation of the scientific-research base at Rzeszów University of Technology" co-financed by the European Union from the European Regional Development Fund within Regional Operational Programme for Podkarpackie Region for the years 2007-2013, I. Competitive and innovative economy, 1.3 Regional innovation system.
\normalsize







\bibliographystyle{plain}   

\bibliography{cissic-a,cissic-b}

\end{document}